\documentclass[11pt,a4paper]{article}% For Computer Modern Fonts, or,
\usepackage{latexsym,amssymb, amsthm}
\usepackage[centertags]{amsmath}
\usepackage{epsfig,pifont}
\usepackage[inactive]{srcltx}  %SRC Specials for DVI search

\usepackage{graphics,helvet,times,mathptm,epic,eepic}

\usepackage{color}

\usepackage{newlfont}

\usepackage{amsfonts,bbm}

\newenvironment{keywords}{ \noindent {\small\bf Key Words}:}{ }

\def\bd{\begin{description}}
\def\ed{\end{description}}

\def\beq{\begin{equation}}
\def\eeq{\end{equation}}
\def\bea{\begin{eqnarray}}
\def\eea{\end{eqnarray}}
\def\beas{\begin{eqnarray*}}
\def\eeas{\end{eqnarray*}}

\theoremstyle{remark}

\begin{document}
%\begin{article}

\title{\textbf{Counting systems and\\
the First Hilbert problem}}

%% Author name
\newcommand{\nms}{\normalsize}
\author{  \bf Yaroslav D. Sergeyev\footnote{Yaroslav D. Sergeyev, Ph.D., D.Sc., is Distinguished
Professor at the University of Calabria, Rende, Italy.
 He is also Full Professor (part-time contract) at the N.I.~Lobatchevsky State University,
  Nizhni Novgorod, Russia  and Affiliated Researcher at the Institute of High Performance
Computing and Networking of the National Research Council of
Italy.
 {\tt yaro@si.deis.unical.it}, \,\, {\tt
http://wwwinfo.deis.unical.it/$\sim$yaro} } \footnote{The author
thanks the anonymous reviewers for their useful suggestions. }}
\date{}

\maketitle

\vspace{1cm}

\begin{abstract}

The First Hilbert problem is studied in this paper by applying two
instruments: a new methodology distinguishing between mathematical
objects and mathematical languages used to describe these objects;
and a new   numeral system allowing one to express different
infinite numbers and to use these numbers for measuring infinite
sets. Several counting systems are taken into consideration. It is
emphasized in the paper that different mathematical languages can
describe mathematical objects (in particular, sets and the number
of their elements) with different accuracies. The traditional and
the new approaches are compared and discussed.

 \end{abstract}

 \vspace{1cm}

\begin{keywords}
The First Hilbert problem;  numeral systems; Pirah\~{a} counting
system;  relativity of mathematical languages; infinite numbers.
 \end{keywords}

% AMS subject classification: 68U, 40A, 03E,11B

%\newpage
 \vspace{1cm}

\section{Introduction}
\label{s1_hilbert}

The First Hilbert problem    belongs to the list of famous 23
mathematical problems that  were announced by David Hilbert during
his lecture delivered before the {I}nternational {C}ongress of
{M}athematicians at {P}aris in 1900 (see \cite{Hilbert}). The
problem has the following formulation.
\begin{quote}
 \textit{\textbf{Cantor's problem of the cardinal number of the continuum }}\\ \\
Two systems, i.e., two assemblages of ordinary real numbers or
points, are said to be (according to Cantor) equivalent or of
equal cardinal number, if they can be brought into a relation to
one another such that to every number of the one assemblage
corresponds one and only one definite number of the other. The
investigations of Cantor on such assemblages of points suggest a
very plausible theorem, which nevertheless, in spite of the most
strenuous efforts, no one has succeeded in proving. This is the
theorem:

Every system of infinitely many real numbers, i.e., every
assemblage of numbers (or points), is either equivalent to the
assemblage of natural integers, $1, 2, 3, \ldots $ or to the
assemblage of all real numbers and therefore to the continuum,
that is, to the points of a line\label{p:hilbert_1}; as regards
equivalence there are, therefore, only two assemblages of numbers,
the countable assemblage and the continuum.

From this theorem it would follow at once that the continuum has
the next cardinal number beyond that of the countable assemblage;
the proof of this theorem would, therefore, form a new bridge
between the countable assemblage and the continuum.

Let me mention another very remarkable statement of Cantor's which
stands in   closest connection with the theorem mentioned and
which, perhaps, offers the key to its proof. Any system of real
numbers is said to be ordered, if for every two numbers of the
system it is determined which one is the earlier and which the
later, and if at the same time this determination is of such a
kind that, if $a$ is before $b$ and $b$ is before $c$, then $a$
always comes before $c$. The natural arrangement of numbers of a
system is defined to be that in which the smaller precedes the
larger. But there are, as is easily seen infinitely many other
ways in which the numbers of a system may be arranged.

If we think of a definite arrangement of numbers and select from
them a particular system of these numbers, a so-called partial
system or assemblage, this partial system will also prove to be
ordered. Now Cantor considers a particular kind of ordered
assemblage which he designates as a well ordered assemblage and
which is characterized in this way, that not only in the
assemblage itself but also in every partial assemblage there
exists a first number. The system of integers $1, 2, 3, \ldots $
in their natural order is evidently a well ordered assemblage. On
the other hand the system of all real numbers, i.e., the continuum
in its natural order, is evidently not well ordered. For, if we
think of the points of a segment of a straight line, with its
initial point excluded, as our partial assemblage, it will have no
first element.

The question now arises whether the totality of all numbers may
not be arranged in another manner so that every partial assemblage
may have a first element, i.e., whether the continuum cannot be
considered as a well ordered assemblage —- a question which Cantor
thinks must be answered in the affirmative. It appears to me most
desirable to obtain a direct proof of this remarkable statement of
Cantor's, perhaps by actually giving an arrangement of numbers
such that in every partial system a first number can be pointed
out.
\end{quote}

Thus, this problem asks about  two questions: (i) Is there a set
whose cardinality is strictly between that of the natural numbers
and that of the real numbers? (ii) Can the continuum of real
numbers be considered a well ordered set?

The question (i) has been stated by Cantor and its answer `no' has
been advanced as the Continuum Hypothesis. G\"odel (see
\cite{Godel}) proved in 1938 that the generalized continuum
hypothesis is consistent relative to Zermelo--Fraenkel set theory.
In 1963, Paul Cohen (see \cite{Cohen}) showed that its negation is
also consistent. These results are not universally accepted as the
final solution to the Continuum Hypothesis  and this area remains
an active topic of contemporary research (see, e.g.,
\cite{Woodin_1,Woodin_2}). The modern point of view on the second
question tells that it is related to Zermelo's axiom of choice
that was demonstrated to be independent of all other axioms in set
theory, so there appears to be no universally valid solution to
the question (ii) either.

In this paper, we look at   the First Hilbert problem using a new
approach  introduced in \cite{Sergeyev}--\cite{Dif_Calculus}. The
point of view on infinite and infinitesimal quantities proposed in
it uses   two methodological ideas borrowed from the modern
Physics: relativity and interrelations   between the object of an
observation and the tool used for this observation. The latter is
applied to  mathematical languages used to describe mathematical
objects and the objects themselves.

In the next section, we give a description of the methodological
principles of the new approach to the interplay between the
counting systems and the objects these systems express (a
comprehensive introduction and examples of its usage can be found
in \cite{www,informatica}).

\section{Counting systems and the relativity of mathematical \\ languages}
\label{s2_hilbert}

The results proved for the First Hilbert problem by G\"odel and
Cohen discuss  how the Continuum Hypothesis is connected to
Zermelo--Fraenkel set theory. Thus,  the problem has been
considered in two logical steps: (i) efforts have been made to
understand what is an accurate definition of the concept `set';
(ii)   the hypothesis was studied using this definition. Such a
way of reasoning is quite common in Mathematics. The idea of
finding an adequate set of axioms for one or another field of
Mathematics is among the most attractive goals for mathematicians.
Usually, when it is necessary to define a concept or an object,
logicians try to introduce a number of axioms describing the
object. However, this way is fraught with danger because of the
following reasons.

First, when we describe a mathematical object or concept we are
limited by the expressive capacity of the language we use to make
this description. A   richer language allows us to say more about
the object and a weaker language -- less. Thus, the incessant
development of the mathematical (and not only mathematical)
languages leads to the continual necessity of a transcription and
the further specification of axiomatic systems. Second, there is
no guarantee that the chosen axiomatic system defines
`sufficiently well' the required concept and a continual
comparison with practice is required in order to check the
goodness of the accepted set of axioms. However, there cannot be
again any guarantee that the new version will be the last and
definitive one because we do not know which new facts related to
the studied object  will be discovered in the future. Finally, the
third limitation has been discovered by G\"odel in his two famous
incompleteness theorems (see~\cite{Godel_1931}).

Thus, if we return to the subject of the First Hilbert problem, we
should observe that people successfully measured finite sets for
centuries  without having a `precise' definition of the notion
`set'. Moreover, the possibility itself of giving such a
definition is questionable (in fact, scientific debates in this
area persist (see, e.g., \cite{Woodin_1,Woodin_2})). In this
paper, we propose a new point of view on the First Hilbert problem
by concentrating our attention on various processes of counting
and mathematical languages used for this purpose.

We start   by noticing that in linguistics the relativity of
descriptions of the world around has been formulated in the form
of the Sapir–-Whorf thesis
(see~\cite{Whorf,Gumperz_Levinson,Hurford,Lucy,Sapir}) also known
as the `linguistic relativity  thesis'. As becomes clear from  its
name, the thesis does not accept the idea of the universality of
language and  challenges the possibility of perfectly representing
the world with language, because it implies that the mechanisms of
any language condition the thoughts of its speakers.

In this paper, we study the relativity   of mathematical languages
in situations where they are used to observe and to  describe
finite and  infinite sets. Let us first illustrate the concept of
the relativity of mathematical languages by the following example.
In his study published in \textit{Science} (see \cite{Gordon}),
Peter Gordon describes a primitive tribe living in Amazonia --
Pirah\~{a} -- that uses a very simple numeral system\footnote{We
remind that \textit{numeral}  is a symbol or group of symbols that
represents a \textit{number}. The difference between numerals and
numbers is the same as the difference between words and the things
they refer to. A \textit{number} is a concept that a
\textit{numeral} expresses. The same number can be represented by
different numerals. For example, the symbols `10', `ten', and `X'
are different numerals, but they all represent the same number
(see \cite{Hurford} for a detailed discussion). } for counting:
one, two, `many'. Sometimes they use the numerals I, II, and III
to indicate these numbers. For Pirah\~{a}, all quantities larger
than two are just `many' and such operations as 2+2 and 2+1 give
the same result, i.e., `many'. By using their weak numeral system
Pirah\~{a} are not able to see, for instance, numbers 3, 4, and 5,
to execute arithmetical operations with them, and, in general, to
say anything about these numbers because in their language there
are neither words nor concepts for that. As a consequence, when
they observe a set having three elements and another set having 4
elements their answer is that both sets have `many' elements.

It is important to emphasize that their answer `many'  is correct
in their language and if one is satisfied with its accuracy, it
can be used (and \textit{is used} by Pirah\~{a}) in practice. Note
that also for us, people knowing that $2+1=3$ and $2+2=4$, the
result of Pirah\~{a} is not wrong, it is just \textit{inaccurate}.
Thus, if one needs a more precise result than `many', it is
necessary to introduce a more powerful mathematical language (a
numeral system in this case) allowing one to express the required
answer in a more accurate way. By using modern powerful numeral
systems where additional numerals  for expressing the numbers
`three' and `four' have been introduced, we can notice that within
`many' there are several objects and the numbers 3 and 4 are among
them. Therefore, these numbers can be used for various purposes,
in particular, for  counting the number of elements in sets.

Among other things  this example   shows us   that a mathematical
language can contain numerals required to formulate a question but
not the numerals that can express the corresponding answer with
the desired accuracy. Thus, the choice of the mathematical
language depends on the practical problem that is to be   solved
and on the accuracy required for such a solution. In dependence of
this accuracy, a numeral system  that would be able to express the
numbers composing the answer (and the intermediate computations,
if any) should be chosen.

Such a situation  is typical for   natural sciences where it is
well known that instruments bound and influence results of
observations. When physicists see a black dot in their microscope
they cannot say: The object of the observation \textit{is} the
black dot. They are obliged to say: the lens used in the
microscope allows us to see the black dot and it is not possible
to say anything more about the nature of the object of the
observation until we  change the instrument - the lens or the
microscope itself - by a more precise one. Then, probably, when a
stronger lens will be put in the microscope, physicists will be
able to see that the object that seemed to be one dot  consists of
two dots.

Note that both results (one dot and two dots) correctly represent
the reality with the accuracy of the chosen instrument of the
observation (lens). Physicists decide the level of the precision
they need and obtain a result depending on the chosen level of the
accuracy. In the moment when they put a lens in the microscope,
they have decided the minimal (and the maximal) size of objects
that they will be able to observe. If they need a more precise or
a more rough answer, they change the lens in their microscope.
Analogously, when mathematicians have decided which mathematical
languages (in particular, which numeral systems) they will use in
their work, they have decided which mathematical objects they will
be able to observe and to describe.

In natural sciences, there always exists the triad -- the
researcher, the object of investigation, and tools used to observe
the object -- and  the instrument used to observe the object
bounds and influences results of observations. The same happens in
Mathematics  studying numbers and objects that can be constructed
by using numbers. Numeral systems used to express numbers are
among  instruments of observations used by mathematicians. The
usage of powerful numeral systems gives the possibility of
obtaining more precise results in Mathematics in the same way as
the  usage of a good microscope gives the possibility to obtain
more precise results in Physics.

Let us now return to Pirah\~{a} again.
%and make the first relevant
%observation with respect to the First Hilbert problem.
Their
numeral system   has another interesting feature particularly
interesting in the context of  the study presented in this paper:
\beq
 \mbox{`many'}+ 1= \mbox{`many'},   \hspace{3mm}
\mbox{`many'} + 2 = \mbox{`many'}, \hspace{3mm} \mbox{`many'}+
\mbox{`many'} = \mbox{`many'}.
 \label{Turing_1}
 \eeq
 These relations are very
familiar to us  in the context of our views on infinity used in
the   Calculus
 \beq
  \infty + 1= \infty,    \hspace{1cm}    \infty
+ 2 = \infty, \hspace{1cm}    \infty + \infty = \infty.
 \label{Turing_2}
 \eeq
Thus, the modern mathematical numeral systems allow us to
distinguish   larger (with respect to Pirah\~{a}) finite numbers
but  when we speak about infinite numbers, they  give results
similar to those obtained by Pirah\~{a}. Formulae (\ref{Turing_1})
and (\ref{Turing_2}) lead us to the following observation:
\textit{Probably our difficulty in working with infinity is not
connected to the nature of infinity but is a result of inadequate
numeral systems used to express infinite numbers.}

This remark is   important with respect to the First Hilbert
problem because its formulation and its investigations executed so
far  used   mathematical instruments developed by Georg Cantor
(see \cite{Cantor}) who has shown that there exist infinite sets
having different number of elements. In particular, the First
Hilbert problem deals with  two kinds of infinite sets --
countable sets and the continuum. Cantor has proved, by using his
famous diagonal argument, that the cardinality, $\aleph_0$, of the
set, $\mathbb{N}$, of natural numbers is less than the
cardinality, $C$, of real numbers    $x \in[0,1]$.

Cantor  has also developed an arithmetic for the infinite cardinal
numbers. Some of the operations of this arithmetic including
numerals  $\aleph_0$ and $C$ are given below:
 \beq
  \aleph_0 + 1
\,\,\,\, =     \aleph_0,  \hspace{1cm} \aleph_0 + 2 \,\,\,\, =
\aleph_0,  \hspace{1cm}
 \label{Hilbert_1}
 \eeq
\beq
 C + 1 \,\,\,\, =     C, \hspace{1cm} C + 2 \,\,\,\, = C,
\hspace{1cm}    C  + \aleph_0 \,\,\,\, = C.
 \label{Hilbert_2}
 \eeq
Again, it is possible to see a clear similarity with the
arithmetic operations used in the numeral system of Pirah\~{a}.
This similarity becomes even more pronounced when one looks at the
numeral system of another Amazonian tribe -- Munduruk\'u  -- who
fail in exact arithmetic with numbers larger than 5 but are able
to compare and add large approximate numbers that are far beyond
their naming range (see \cite{Pica}). Particularly, they use the
words `some, not many' and `many, really many' to distinguish two
types of large numbers. It is sufficient to substitute in
(\ref{Hilbert_1}), (\ref{Hilbert_2}) these two words with the
numerals $\aleph_0$ and $C$, respectively, to have an idea about
their arithmetic.

Thus, the point of view on Mathematics, and, in particular, on
numbers and sets of   numbers ($\mathbb{N}, \mathbb{Z},
\mathbb{R},$ etc.) used in this paper follows natural sciences and
consists of the following. There exist mathematical objects (e.g.,
numbers and sets of numbers) that are   objects of the
observation. There exist different numeral systems (instruments of
the observation) allowing us to observe certain numbers through
numerals belonging to these numeral systems, to execute certain
operations with them, and, in particular, to use these numerals
for measuring sets of numbers. The researcher (the observer) is
able to view only those numbers that are visible through numerals
available in one or another numeral system. Sometimes a fixed
numeral system, $\mathcal{S}_1$, allows the researcher to obtain a
precise answer that cannot be improved by the usage of a more
sophisticated numeral system, $\mathcal{S}_2$, and sometimes the
introduction of $\mathcal{S}_2$ allows him/her to improve the
accuracy of the result.

For example, let us indicate the numeral system of Pirah\~{a}  as
$\mathcal{S_P}$ and the numeral system of Munduruk\'u as
$\mathcal{S_M}$. Then the accuracy of expression of the numbers 1
and 2 by $\mathcal{S_P}$ cannot be improved by switching to the
usage of $\mathcal{S_M}$. However, the accuracy of the answer
`many' can be improved by $\mathcal{S_M}$ when numbers less than
six are observed. The system $\mathcal{S_M}$, in its turn, starts
to give inaccurate answers when one starts to use it to observe
numbers larger than five. Modern numeral system (including
numerals $\infty$, $\aleph_0$, $\aleph_1$,  $C$, $\varpi$, etc.)
can improve this accuracy for large finite numbers but start to
give inaccurate answers when infinite (and infinitesimal) numbers
should be observed.

This clear separation between mathematical objects and
mathematical lan\-gua\-ges used to observe and to describe these
objects means, in particular,   that from this point of view,
\textit{axiomatic systems do not define mathematical objects but
just determine formal rules for operating with certain numerals
reflecting   only some (not all)  properties of the studied
mathematical objects}. For example, axioms for real numbers are
considered together with a particular numeral system $\mathcal{S}$
used to write down numerals and are viewed as practical rules
(associative and commutative properties of multiplication and
addition, distributive property of multiplication over addition,
etc.) describing operations with the numerals. The completeness
property is interpreted as a possibility to extend $\mathcal{S}$
with additional symbols (e.g., $e$, $\pi$, $\sqrt{2}$, etc.)
taking care of the fact that the results of computations with
these symbols  agree with the facts observed in practice. As in
Physics, it becomes impossible to say anything about the object of
the observation without the usage of an instrument of the
observation. In case of mathematical objects,  the used
mathematical languages   (including numeral systems)  are among
these instruments.

Another important consequence of the understanding  of the
existence of mathematical objects and their separation from tools
of the observation consists of the fact that it is necessary to be
very careful when one speaks about sets of numbers. In particular,
when we speak about sets (finite or infinite) we should take care
about tools used to describe a set. In order to introduce a set,
it is necessary to have a language (e.g., a numeral system)
allowing us to describe its elements and the number of the
elements in the set. For instance, the set
 \beq
  A=\{1, 2, 3, 4, 5\}
\label{4.1.deriva_0}
 \eeq
 cannot be defined using the mathematical
language of Pirah\~{a}.

Such words as `the set of all finite numbers' do not define a set
completely from our point of view, as well. It is always necessary
to specify which instruments are used to describe (and to observe)
the required set and, as a consequence, to speak about `the set of
all finite numbers expressible in a fixed numeral system'. For
instance, for Pirah\~{a} `the set of all finite numbers'  is the
set $\{1, 2 \}$ and for   Munduruk\'u   `the set of all finite
numbers' is the set $A$ from (\ref{4.1.deriva_0}). In modern
numeral systems there exist similar limitations and hence our
possibilities to write down numbers are limited (e.g., it is not
possible to write down  in a positional numeral system any number
that has $10^{100}$ digits)\footnote{It is interesting that this
problem was clear to Archimedes. In his work \textit{The Sand
Reckoner}, in order to express huge numbers he had invented new
numerals. It was necessary because the ancient Greeks used a
simple system for writing numbers using 27 different letters   of
the Greek alphabet. This system did not allow Archimedes to write
down huge numbers he wished to work with.}. We emphasize again
that, as it happens in Physics, the instrument used for an
observation bounds the possibility of the observation. It is not
possible to say how we shall see an object of our observation if
we have previously not clarified which instruments are to be used
to execute the observation.

Thus, we have emphasized in this section that there exists the
crucial difference between numbers  and numeral systems used to
observe them and, analogously, between sets of numbers and sets of
numerals used to observe them. Advanced contemporary numeral
systems enable us to distinguish within `many'   various large
finite numbers. As a result, we are able to use large finite
numbers in computations, to construct mathematical models
involving them, and, in particular, to measure large finite sets.
Analogously, it can be expected that if we would be able to
distinguish more infinite numbers we could execute computations
with them and measure infinite sets.

\section{Mathematical objects
discussed in the First Hilbert \\ problem and the accuracy of
mathematical languages used for their description}
\label{s3_hilbert}

Let us  start by making some observations with respect to the
formulation of the part (i) of the First Hilbert problem (see
Introduction).  It uses the mathematical language developed by
Cantor and considers real numbers, natural numbers, cardinal
numbers, countable sets, and the continuum. In the formulation
there is no any distinction between objects of the observation
(infinite sets and numbers) and the mathematical language   used
for the observation (in particular, the  words `countable sets'
and `continuum'). Let us reconsider mathematical objects and
numeral systems involved in the formulation of the problem by
concentrating our attention on this distinction.

The real numbers are understood in the formulation  (see
page~\pageref{p:hilbert_1}) to be `the points of a line'. Then,
the existence of the continuum has been  demonstrated by Cantor
using the diagonal argument working with a positional
representation of points
 \beq
(.a_{1} a_{2} a_{3} \ldots a_{i-1}a_{i}a_{i+1} \ldots )_b
 \label{Hilbert_3}
 \eeq
in the unitary interval where the symbol $b$ indicates a finite
radix of the positional system and symbols $a_{i}, 0 \le a_{i} \le
b-1,$ are called digits. Thus, the object of the observation --
the set of real numbers over the unitary interval -- is observed
through the positional numeral system (\ref{Hilbert_3}) with the
radix $b$.

Then the numeral system of infinite cardinals is used to observe
infinite numbers and to measure infinite sets. In particular, two
cardinal infinite numbers expressed by the numerals $\aleph_0$ and
$C$ are   used to measure  the infinite sets of natural and real
numbers in the formulation of the First Hilbert problem.

Let us make some important observations. First,  Cantor has
proved, that the number of elements of the set $\mathbb{N}$ is
smaller  than the number of real numbers at the interval $[0,1)$
\textit{without calculating the latter}. Second, he worked with
the numerals (\ref{Hilbert_3}) supposing that they represent
\textit{all real numbers} over the unitary interval. Third, he has
also proved   that the introduced numeral $C$ is the cardinality
also of the whole real line.

However, as it has been discussed in the previous section, the
arithmetic using numerals $\aleph_0$ and $C$ clearly indicates
that the accuracy of this numeral system is not satisfactory when
it is necessary to measure infinite sets (e.g.,    addition of an
element to any infinite  set cannot be registered, see
(\ref{Hilbert_1}), (\ref{Hilbert_2})). This numeral system allows
us to distinguish sets having `some, not many' elements from the
sets having `many, really many' elements but does not give any
possibility to distinguish the measure of sets $D$ and $\tilde{D}$
if both $D$ and $\tilde{D}$ have `some, not many' elements or both
$D$ and $\tilde{D}$ have `many, really many' elements.

 Let us
consider, for example,   the set of integers, $\mathbb{Z}$, and
the set of natural numbers, $\mathbb{N}$. The result saying that
both of them have the same cardinal number $\aleph_0$ can be
viewed as an indication of the fact that the numeral system of
cardinals is too weak to discern the difference of the number of
elements of these two sets. An analogous comment can be made with
respect to the fact that $C$ is the cardinality of both the set of
all   real numbers, $\mathbb{R}$, and the   set of real numbers
over the unitary interval being a proper subset of $\mathbb{R}$.

Thus, a more powerful numeral system that would allow us to
capture these differences is required to measure infinite sets.
Recently, a new   numeral system has been developed to express
finite, infinite, and infinitesimal numbers in a unique framework
(a rather comprehensive description of the new methodology  can be
found  in the survey \cite{informatica}). The main idea consists
of measuring infinite and infinitesimal quantities   by different
(infinite, finite, and infinitesimal) units of measure using a new
numeral.

A new infinite unit of measure   has been introduced  for this
purpose as the number of elements of the set $\mathbb{N}$ of
natural numbers. It is expressed by a new numeral \ding{172}
called \textit{grossone}. The    new approach   does not
contradict Cantor and can be viewed as an evolution of his deep
ideas regarding the existence of different infinite numbers in a
more applied and precise way.

In the new numeral system using \ding{172}, it becomes possible to
express a variety of numerals   that allow one to observe infinite
positive integers and to use them for measuring
 certain infinite sets. As a consequence, the new
numeral system allows us to improve the accuracy of measuring
infinite sets in comparison with the numeral system using
cardinals introduced by Cantor. For instance, within the countable
sets it becomes possible to distinguish the following sets having
the different number  of elements: the set of even numbers has
$\frac{\mbox{\ding{172}}}{2}$ elements; the set of integers,
$\mathbb{Z}$, has $2\mbox{\ding{172}\small{+}}1$ elements;  the
set $\mathbb{N}_{-} =\mathbb{N}\setminus \{x\},$   $x
\in\mathbb{N},$ has $\mbox{\ding{172}}-1$ elements;   the set
$\mathbb{N}_{+} =\mathbb{N} \bigcup \{y\},$   $y  \notin
\mathbb{N},$ has $\mbox{\ding{172}}+1$ elements; and the set
 \[
 P  =
\{
  (a_1, a_2) : a_1 \in   \mathbb{N}, \,\, a_2 \in   \mathbb{N}
 \}
  \]
  has $\mbox{\ding{172}}^2$ elements.  Another example of
the usage of the new numerals can be given for the set $B = \{
3,4,5, 69 \} \cup (B_1 \cap B_2 )$, where $B_1 \subset
\mathbb{N},\,\,   B_2 \subset \mathbb{N}$ and
\[
B_1 = \{ 4, 9, 14, 19, 24,    \ldots\}, \hspace{1cm}
 B_2 = \{ 3,
14, 25, 36, 47,  \ldots\}.
\]
As it has been shown in \cite{informatica}, the set $B$ has
$\frac{\mbox{\ding{172}}}{55}+3$ elements.

With respect to the sets that were indicated by Cantor as having
cardinality of the continuum the new approach also allows one to
obtain more precise results and to notice that among sets having
cardinality $C$   there exist   sets having the different infinite
number of elements. First of all, it can be shown (see
\cite{informatica}) that numerals (\ref{Hilbert_3}) do not
represent \textit{all} real numbers within the unitary interval.
They represent only those numbers that can be expressed by the
numerals (\ref{Hilbert_3}) that   have no more than grossone
digits. For instance, the irrational number $\frac{\sqrt {e}}{4}$
cannot be expressed in it because even grossone digits are not
sufficient to express $e$ in the positional numeral system (see
\cite{informatica}).

Moreover, the new approach allows us to show that the numbers,
$n_b$, of all numerals expressible    in the form
(\ref{Hilbert_3}) are different depending on the radix $b$ and to
calculate that $n_b=b^{\mbox{\tiny{\ding{172}}}}$. For instance,
the binary numeral system   ($b=2$) contains
$n_2=2^{\mbox{\tiny{\ding{172}}}}$ numerals. This is smaller than,
for example, the number of numerals expressible in the decimal
positional system having $n_{10}=10^{\mbox{\tiny{\ding{172}}}}$
numerals. As a result, the binary numeral system can express less
numbers than the decimal one. It is important to emphasize that
for any finite $b>1$ all of the obtained numbers,
$b^{\mbox{\tiny{\ding{172}}}}$, are larger than \ding{172} being
the number of elements of the set, $\mathbb{N}$, of natural
numbers. This fact can be viewed as a further specification of the
results obtained by Cantor.

We can also consider positional numeral systems having an infinite
number of positions $k$ in (\ref{Hilbert_3}) such that
 $k<\mbox{\ding{172}}$. Then these sets will have less than
$b^{\mbox{\tiny{\ding{172}}}}$ numerals. For instance, for
$k=\mbox{\ding{172}}/2$, we have the following numeral system
 \[
(.a_{1} a_{2} a_{3} \ldots a_{0.5
\mbox{\tiny{\ding{172}}}-1}a_{0.5\mbox{\tiny{\ding{172}}}} )_b
  \]
that contains $b^{\mbox{\tiny{\ding{172}}}/2}$ numerals.

Another important particular case can be obtained from the
solution to the equation $b^{x}=\mbox{\ding{172}}$ that tells us
that the positional numeral systems having $k_1=[\log_b
(\mbox{\ding{172}})]$ and $k_2=[\log_b (\mbox{\ding{172}})]+1$
positions in their numerals can express $y_1 \le
\mbox{\ding{172}}$, and $y_2 > \mbox{\ding{172}}$ numbers,
respectively, where $[u]$ is the integer part of $u$. The
positional numeral system with the numerals having $k_1$ digits is
of special interest in the context of the Continuum Hypothesis and
we call it \textit{critical value}. This terminology is used
because even though in this case we deal with a positional system
with numerals consisting of an infinite number of digits, the
quantity, $y_1$, of numerals expressible in this numeral system is
less than the number of elements of the set of natural numbers. By
a complete analogy critical values can be calculated not only with
respect to grossone but for other infinite integers, as well. For
instance, this can be done for $\mbox{\ding{172}}/2$ (the number
of elements of the set of even numbers), for
$2\mbox{\ding{172}}+1$ (the number of elements of the set of
integers), etc.

Let us notice now that other numeral systems can be used to
express points of an interval (see \cite{informatica}) and that
the new methodology allows us to introduce a more  physical
concept of continuity (see \cite{Dif_Calculus}). Recall that in
Physics, the `continuity' of an object is relative. For example,
if we observe a table by eye, then we see it as being continuous.
If we use a microscope for our observation, we see that the table
is discrete. This means that \textit{we decide} how to see the
object, as a continuous or as a discrete, by the choice of the
instrument for the observation. A weak instrument -- our eyes --
is not able to distinguish its internal small separate parts
(e.g., molecules) and we see the table as a continuous object. A
sufficiently strong microscope allows us to see the separate parts
and the table becomes discrete but each small part now is viewed
as continuous.

In contrast, in   traditional mathematics, any mathematical object
is either continuous or discrete. For example, the same function
cannot be   both  continuous and discrete. Thus, this
contraposition of discrete and continuous in the traditional
mathematics does not reflect  properly the physical situation that
we observe in practice. The numeral system including grossone
gives us the possibility to develop a new theory of continuity
that is closer to the physical world and better reflects the new
discoveries made by the physicists (see \cite{Dif_Calculus}).

Let us pass now from the unitary interval to the whole real line.
To observe real numbers (the object of the observation) it is
necessary first to choose a numeral system (the instrument of the
observation) to represent them. After this choice has been done,
it becomes possible to count the number of numerals expressible in
the fixed numeral system and to understand how many real numbers
can be expressed in this numeral system. For instance, let us
consider numerals
 \beq
 \pm(  \ldots a_2 a_1 a_0 .
 a_{1} a_{2} a_{3} \ldots   )_b
   \label{Hilbert_4}
 \eeq
that are expressed in the positional system with the radix   $b$.
Their number is equal to $2b^{2\mbox{\tiny{\ding{172}}}}$ if we
want to use in (\ref{Hilbert_4}) sequences of integer and
fractional digits consisting of \ding{172} elements each (see
\cite{informatica}). Note that $2b^{2\mbox{\tiny{\ding{172}}}}$ is
the number of different numerals. In this numeral system the
number zero can be represented by two different numerals
\[
 -\underbrace{000\ldots000}_{\mbox{\tiny{\ding{172} digits}}}.
 \underbrace{000\ldots000}_{\mbox{\tiny{\ding{172} digits}}},
\hspace{1cm}+\underbrace{000\ldots000}_{\mbox{\tiny{\ding{172}
digits}}}.\underbrace{000\ldots000}_{\mbox{\tiny{\ding{172}
digits}}}.
\]

Suppose not that we want to change the numeral system for
representation of real numbers and decide to work, e.g., with the
numerals
 \beq
 \pm(.
 a_{ 1} a_{ 2} a_3 \ldots    )_b \cdot b^{\pm( p_{ 1} p_{ 2} p_3  \ldots   )_b}
 \label{Hilbert_5}
 \eeq
where   digits $0 \le a_i < b$ and $0 \le p_i < b$, and both
sequences $\{ a_{ i}\}$ and $\{ p_{ i}\}$ consist  of \ding{172}
elements each. Then it is easy to show that  this numeral system
contains $4b^{2\mbox{\tiny{\ding{172}}}}$ different numerals.

 We emphasize again  that any  numeral system
allows us to express precisely the number of elements only of
certain sets. For the set (\ref{4.1.deriva_0}) the numeral system,
$\mathcal{S_M}$, of Munduruk\'u allows us to do this. For the sets
$\mathbb{N}$ and $\mathbb{Z}$ the numeral system including
grossone is sufficient. To express the number of elements of the
set of real numbers even this numeral system is not sufficient.
Even though it   allows us to specify certain infinite subsets of
the set of real numbers and to count their number of elements, it
is too weak for both expressing the number of elements of
$\mathbb{R}$ and expressing irrational numbers. Note that with
respect to the latter it behaves as  all other existing numeral
systems which use one of the following two alternatives: (a)   use
approximations of irrational numbers expressible in the chosen
numeral system;
  (b)   use numerals created ad hoc for the required
irrational numbers, e.g. $e, \pi, \sqrt{2}, \sqrt[5]{3}$, etc.

Let us consider now the second part of the First Hilbert Problem
and give some comments upon. Again, from our point of view, its
difficulty  consists of the fact that in its formulation there is
no   distinction between the object of the observation and the
numeral systems used for the observation.

In our analysis, we first  emphasize that in order to compare two
real numbers it is necessary to have a numeral system allowing us
to express both of them. Then again, as it was in Physics with the
choice of the lens for the microscope, we have decided which
numbers can be compared at the moment when we have chosen our
numeral system.
  This   is very important
from the methodological point of view: we are not able to compare
numbers if we have no numerals allowing us to express these
numbers (at least in an approximative way).

For instance, let us consider the set, $\overline{\mathbb{N}}$, of
all natural numbers   that we are able to write down using all
known numeral systems. Then, we are not able to write down, to
compare, and to order numbers belonging to the set $\mathbb{N}
\setminus \overline{\mathbb{N}}$. This set can be defined   but we
are not able to indicate any of its elements, including the first
one. Thus, it is possible to speak about the existence of this set
but it is impossible to execute operations with its elements
because we have no tools to express them.

Note also the importance of the fact that if we want to compare
two numbers accurately, they should be expressible in the same
numeral system. If they are expressed by two different numeral
systems then the problem of a transcription from one to another
system arises and it is possible that such a precise transcription
is impossible (e.g., it is shown in \cite{informatica} that $e$
cannot be written in any positional system with a finite radix $b$
even if one has an infinite number of digits in this positional
record). In such cases only an approximate comparison can be done
($e \approx 2.7$).

Let us now return to the formulation of the problem where Hilbert
writes `On the other hand the system of all real numbers, i.e.,
the continuum in its natural order, is evidently not well ordered.
For, if we think of the points of a segment of a straight line,
with its initial point excluded, as our partial assemblage, it
will have no first element.' In these phrases Hilbert speaks about
numbers expressible by numerals (\ref{Hilbert_3})  as about all
real numbers. This is not the case as it has been shown in the
analysis above. Several numeral systems that can be used to
express real numbers have been studied in \cite{informatica} in
detail.

With respect to the second phrase of Hilbert it is possible to
make the following comment. Without loss of generality, we
consider the unitary interval as   a segment of a straight line
and we use numerals (\ref{Hilbert_3})   to express numbers Hilbert
speaks about. It has been shown in \cite{informatica} for these
numerals (other numeral systems have been   studied in
\cite{informatica}, as well) that all of them represent different
numbers. Thus, if we have decided to observe the interval $[0,1]$
by   numerals (\ref{Hilbert_3}) with $b=10$ then only the
following $10^{\mbox{\tiny{\ding{172}}}}$ numbers can be
distinguished within this interval
\[
0.\underbrace{000\ldots000}_{\mbox{\tiny{\ding{172} digits}}},\,\,
0.\underbrace{000\ldots0001}_{\mbox{\tiny{\ding{172} digits}}},
\hspace{3mm}\ldots
\hspace{3mm}0.\underbrace{999\ldots998}_{\mbox{\tiny{\ding{172}
digits}}},\,\, 0.\underbrace{999\ldots999}_{\mbox{\tiny{\ding{172}
digits}}}.
\]
Then, if we decide to exclude from the interval $[0,1]$ the number
zero represented by the numeral
$0.\underbrace{000\ldots000}_{\mbox{\tiny{\ding{172} digits}}}\,$
then the next numeral,
$0.\underbrace{000\ldots0001}_{\mbox{\tiny{\ding{172} digits}}}$,
gives us the first element in the set of remaining numerals (that
has now $b^{\mbox{\tiny{\ding{172}}}}-1$ elements). Obviously, if
the accuracy of the numeral system (\ref{Hilbert_3}) is not
sufficient for a problem we want to solve over the interval
$[0,1]$, then we can decide to add some new numerals to deal with.
For instance, we can decide to add the numeral
$1.\underbrace{000\ldots000}_{\mbox{\tiny{\ding{172} digits}}}\,$
in order to represent the number one exactly. Then,   this numeral
and numerals (\ref{Hilbert_3}) give us the possibility to
distinguish within the interval $b^{\mbox{\tiny{\ding{172}}}}+1$
different numbers.

Therefore, the numbers that can be expressed by the numerals
(\ref{Hilbert_3}) and that we are able to write down can be
ordered. This note is important because we are not able to write
down all numerals of the type (\ref{Hilbert_3}), e.g., the ones
having an infinite number of different digits. However, the first
elements Hilbert speaks about in the text quoted above can be
expressed and ordered. Without loss of generality we give an
example using   the decimal positional system and write down the
first four numerals from this numeral system:
\[
0.\underbrace{000\ldots000}_{\mbox{\tiny{\ding{172} digits}}} <
0.\underbrace{000\ldots0001}_{\mbox{\tiny{\ding{172} digits}}} <
0.\underbrace{000\ldots0002}_{\mbox{\tiny{\ding{172} digits}}} <
0.\underbrace{000\ldots0003}_{\mbox{\tiny{\ding{172} digits}}}.
\]
 By a complete analogy, the  smallest numbers expressible in the
 numeral systems (\ref{Hilbert_4}) and (\ref{Hilbert_5}) can be
 indicated.

\section{Conclusion}
\label{s4_hilbert}

Two new  mathematical instruments have been applied in the paper
to study the First Hilbert problem: the methodology distinguishing
 mathematical objects and mathematical languages used to
describe these objects; and the new  numeral system allowing one
to express different infinite numbers and to use these numbers for
measuring infinite sets.

With respect to the Continuum Hypothesis it has been established
that the mathematical language used by  Cantor and Hilbert to
formulate the problem has two peculiarities. First, it does not
take into consideration the difference between mathematical
objects under the observation and the mathematical language  used
to describe these objects. Second, this language does not allow
one to discern among different numerable sets those sets that have
the different infinite number of elements. The same happens with
respect to different sets having the cardinality of the continuum.

The new numeral system has allowed us to measure certain infinite
sets with a high  accuracy and to express explicitly the number of
their elements. In particular,  the sets mentioned in the
Continuum Hypothesis have been studied in detail and   a
constructive answer to the Hypothesis has been given in the
framework of the new methodology.

It is obligatory to emphasize that both approaches, that of Cantor
and Hilbert and the new one, do not contradict one another. Both
of them   represent the reality but they do it with the different
accuracy. Such a situation  is typical for   natural sciences
where it is well known that instruments bound and influence
results of observations. As a result, the moment when  a
researcher chooses  an instrument for the required observation is
the moment when the accuracy of the observation is determined. We
can illustrate this situation by using again the analogy with the
microscope.

Suppose that with a weak lens a physicist  sees two  dots in the
microscope and with a stronger lens instead of the two dots he/she
sees 45 smaller dots. Both answers are correct but they describe
the object of the observation with the different   accuracy that
is determined by the instrument (the lens) chosen for the
observation. Analogously, when mathematicians have decided which
mathematical languages (in particular, which numeral systems) they
would use in their work, they have decided which mathematical
objects they would be able to observe and to describe.

The language used by   Cantor and Hilbert allows an observer to
see  in the mathematical microscope two dots (numerable sets and
the continuum) and the accuracy of the used lens (the numeral
system using cardinals) does not allow the observer to capture the
presence of   eventual other dots. In addition, the observer does
not realize that he/she has put on the sample stage of the
microscope instead of the entire object (real numbers) just a part
of it (namely, numbers expressible in a positional numeral system
with the infinite number of digits). Moreover, the observer is not
conscious that even the taken part is not always the same because
it  is changed when the radix of the positional numeral system is
changed.

In the   case of the new mathematical language, the observer
distinguishes the lens from the objects of the observation
(numeral systems used to describe sets from the sets themselves)
and has a stronger lens (the numeral system using gross\-one)
allowing him/her to observe instead of two dots   many different
dots (various infinite sets having the different number of
elements). In particular, it becomes possible to distinguish
various sets of numerals expressible by positional systems with
different radixes and different infinite lengths of digits. It has
been shown that for a fixed   radix  a critical length can been
established such that sets of numerals with a length superior to
the critical one have more elements than the set of natural
numbers. For infinite lengths inferior to the critical value the
number of numerals expressible by this positional system  is
inferior to the number of elements of the set of natural numbers.
Analogous critical values can be established for other numerable
sets.

With respect to the second part of the problem dealing with well
ordered sets it has been commented that it is necessary to operate
very carefully with  propositions dealing with the existence of
non trivial mathematical objects and with executing operations
with them. In particular,   if such a proposition supposes the
execution of an operation, then it is necessary to verify the
existence of mathematical tools (e.g.,   numeral systems) allowing
us to express mathematical objects required for this operation
including both the operands and the intermediate and final
results.

It has been emphasized that we are not able to compare numbers if
we have no numerals allowing us to express these numbers (at least
in an approximative way). As a consequence, only those numbers can
be ordered for which  numeral systems allowing one to express them
are known. Finally, it has been shown that minimal numbers
expressible in   positional numeral systems with the different
infinite quantity of digits  can be ordered.

\bibliographystyle{plain}
\bibliography{XBib_First}

\end{document}